\documentclass[11pt, leqno]{article}
\usepackage{amsmath, amsfonts}
\begin{document}

\author{Ajai Choudhry}
\title{Quadratic diophantine equations\\ with applications to quartic equations}
\date{}
\maketitle

\begin{abstract} In this paper we first show that, under certain conditions, the solution of  a single quadratic diophantine equation in four variables $Q(x_1,\,x_2,\,x_3,\,x_4)=0$ can be expressed in terms of bilinear forms in four parameters. We use this result to establish a necessary, though not sufficient, condition for the solvability of the simultaneous quadratic diophantine equations  $Q_j(x_1,\,x_2,\,x_3,\,x_4)=0,\;j=1,\,2,$ and give a method of obtaining their complete solution. In general, when these two equations have a rational solution, they represent an elliptic curve but we show that there are several cases in which their  complete solution may be expressed by  a finite number of parametric solutions and/ or a finite number of primitive integer solutions. Finally we relate the  solutions of the quartic equation $y^2=t^4+a_1t^3+a_2t^2+a_3t+a_4$ to the solutions of a pair of quadratic diophantine equations, and thereby obtain new formulae for deriving rational solutions of the aforementioned  quartic equation starting from one or two known solutions.

\end{abstract}
 
Keywords: bilinear solutions of quadratic diophantine equations; quartic diophantine equation; quartic model of elliptic curve; quartic function made a perfect square.

Mathematics Subject Classification 2010: 11D09, 11D25.

\section{Introduction}
\setcounter{section}{1}

\hspace{0.25in} In  this paper  we study the  solutions of the single quadratic diophantine equation, 
\begin{equation}
Q(x_1,\,x_2,\,x_3,\,x_4)=0, \label{qdgen}
\end{equation}
as well as solutions of a pair of simultaneous quadratic diophantine equations,
\begin{equation}
Q_j(x_1,\,x_2,\,x_3,\,x_4)=0,\;j=1,\,2, \label{qdpair}
\end{equation} 
where $Q(x_1,\,x_2,\,x_3,\,x_4)$ and $Q_j(x_1,\,x_2,\,x_3,\,x_4), \;j=1,\,2,$  are quadratic forms with rational coefficients in four variables $x_i,\,i=1,\,2,\,3,\,4.$ We also  consider the application of these solutions to the quartic diophantine equation,
\begin{equation}
y^2=t^4+a_1t^3+a_2t^2+a_3t+a_4, \label{quartic0}
\end{equation}
where the coefficients $a_i,\,i=1,\,2,\,3,\,4,$ are rational numbers. 

The usual  method of solving equation (\ref{qdgen}), as described in \cite[p. 432]{Di}, consists in first finding a single numerical solution, and then using it to obtain the complete solution. This complete solution is obtained in terms of polynomials of degree two in three parameters. There are very few  equations 
 of type (\ref{qdgen}) for which  the existing literature  gives the complete solution in terms of bilinear forms in four parameters (see, for instance, \cite[p. 15]{Mo}). We show that the complete solution of the quadratic diophantine equation (\ref{qdgen}) can be expressed, under certain conditions,  in terms of bilinear forms in four variables. For instance, the complete solution of the diophantine equation 
\begin{equation}
270x_1^2+76x_1x_2+152x_1x_3-16x_2^2-48x_2x_3-35x_3^2+3x_4^2=0, \label{ex2}
\end{equation}
may be expressed as,
\begin{equation}
\begin{array}{rclrcl}
x_1 &=& 4pm-2qn, &x_2 &=& (117p+4q)m+(2p-65q)n,\\ x_3 &=& -(72p+2q)m-(p-40q)n, &x_4 &=&2qm -pn, \label{solex2}
\end{array}
\end{equation}
where $p,\,q,\,m,\, n,$ are arbitrary parameters.

In fact, in Section 2 we obtain a necessary and sufficient condition for the solution of an equation of type (\ref{qdgen}) to be expressed in terms of bilinear forms and  give a method of obtaining such solutions. 

In Section 3 we obtain a necessary condition for the solvability of the simultaneous equations (\ref{qdpair}) and give a method of obtaining their complete solution when this condition is satisfied. It is well-known that, in general, when the two equations have a rational solution, they represent an elliptic curve. We show, however, that there are numerous cases when the complete solution of these  equations  may be given by one or more parametric solutions and they may have, in addition, a finite number of primitive solutions. As an example, the complete solution of the simultaneous diophantine equations,
\begin{equation}
\begin{array}{r}
  x_1^2-9x_2^2-x_3^2+4x_4^2=0, \\
3x_1^2-30x_1x_2-4x_1x_3-9x_2^2 -12x_2x_3-7x_3^2+12x_3x_4+4x_4^2=0, \label{qdpairex50}
   \end{array} \end{equation}
 is given by 
	two parametric solutions, a linear solution namely,
   \begin{equation}
	x_1 =2m ,\; x_2 =2n,\; x_3 = -2m, \;x_4 =-3n, \label{qdpairex5sol1}
	\end{equation}
   and a solution of degree 3,
   \begin{equation}\begin{array}{lcl}x_1 =6m^2n-108mn^2+270n^3 ,& x_2 = 2m^3-8m^2n+18mn^2-36n^3 , \\
	x_3 =30m^2n-36mn^2+270n^3 ,& x_4 =3m^3-12m^2n+63mn^2+54n^3, \label{qdpairex5sol2}
   \end{array}\end{equation}
   where $m,\,n$ are arbitrary parameters.

In Section 4 we derive formulae giving new solutions of the quartic equation (\ref{quartic0}) starting from one or two known solutions. While there is a vast amount of literature dealing with equation (\ref{quartic0}), the formulae given in Section 4 do not seem to have been obtained earlier. The solutions of (\ref{quartic0}) obtained by these formulae are much simpler than those obtained by existing methods. 

With respect to the solutions of homogeneous diophantine equations of type (\ref{qdgen}) or (\ref{qdpair}), we note that the existence of any rational solution $(x_1,\,x_2,\,x_3,\,x_4)=(\alpha_1,\,\alpha_2,\,\alpha_3,\,\alpha_4)$ implies that $(x_1,\,x_2,\,x_3,\,x_4)=(k\alpha_1,k\,\alpha_2,\,k\alpha_3,\,k\alpha_4), \; k \in {\mathbb Q}\setminus\{0\}$ is also a rational solution, and all such solutions will be considered equivalent. This equivalence class of solutions will be considered as a single solution of the equation(s) under consideration.

\section{The  diophantine equation $Q(x_1,\,x_2,\,x_3,\,x_4)=0$}

\hspace{0.25in}In this section  we first prove a preliminary lemma regarding the diophantine equation
\begin{equation}
 a_1x_1^2+a_2x_2^2+a_3x_3^2+a_4x_4^2=0, 
\end{equation}
and then prove  a theorem that gives a necessary and sufficient condition for the complete solution of this equation  to be expressed in terms of bilinear forms in four parameters. We also obtain this  bilinear solution explicitly  when this condition is satisfied. Next  we prove two   theorems that give conditions for the existence of  bilinear solutions of  the general quadratic equation (\ref{qdgen}) and show how such  solutions can be obtained.  We also give a couple of examples of bilinear solutions of quadratic equations.

\subsection{Bilinear solutions of the equation $Q(x_1,\,x_2,\,x_3,\,x_4)=0$}

\noindent{\bf Lemma 1:} 
 If $a_i,\,i=1,\,2,\,3,\,4,$ are nonzero rational numbers such that $a_1a_2a_3a_4$ is a perfect square, the diophantine equation
\begin{equation}
a_1x_1^2+a_2x_2^2+a_3x_3^2+a_4x_4^2=0 \label{qd4var}
\end{equation}
has a nontrivial solution in integers if and only if 
the diophantine equation
\begin{equation}
a_1y_1^2+a_2y_2^2+a_3y_3^2=0 \label{qd3var}
\end{equation}
has a nontrivial solution in integers. 

\noindent{\bf Proof:} If $a_1a_2a_3a_4=k^2, \; k \in {\mathbb Q}\setminus\{0\}, $ and equation (\ref{qd4var}) has a nontrivial solution in integers,  two of the coefficients $a_i$ must be positive and two must be negative. We assume without loss of generality that $a_1$ and $a_2$ are positive and $a_3$ and $a_4$ are negative. If $x_i=\alpha_i,\,i=1,\,2,\,3,\,4,$ is a nontrivial solution of (\ref{qd4var}), we   write
\begin{equation} \begin{array}{rcl}
\beta_1 &=& a_2(k\alpha_2\alpha_4+a_1a_3\alpha_1\alpha_3),\\
\beta_2 &=& a_1(k\alpha_1\alpha_4-a_2a_3\alpha_2\alpha_3),\\ 
\beta_3 &=& a_1a_2(a_3\alpha_3^2+a_4\alpha_4^2). \label{valbeta}
\end{array}
\end{equation} 
We note that $\beta_3 \neq 0$ since the relation $\beta_3=0$ would imply that both $\alpha_3,\,\alpha_4$ are 0, and hence both $\alpha_1,\,\alpha_2$ must also be zero contradicting the fact that $x_i=\alpha_i$ is a nontrivial solution of (\ref{qd4var}). Now, on using the relation $a_1a_2a_3a_4=k^2$, we get 
\[ 
a_1\beta_1^2+a_2\beta_2^2+a_3\beta_3^2=a_1^2a_2^2a_3(a_3\alpha_3^2+a_4\alpha_4^2)(a_1\alpha_1^2+a_2\alpha_2^2+a_3\alpha_3^2+a_4\alpha_4^2)=0. \]
It follows that  $y_i=\beta_i,\,i=1,\,2,\,3,$ is a nontrivial solution of (\ref{qd3var}).  Conversely, if  
$y_i=\beta_i,\,i=1,\,2,\,3,$ is a given nontrivial solution of (\ref{qd3var}), a nontrivial solution of (\ref{qd4var})
is given by $(x_1,\,x_2,\,x_3,\,x_4)=(\beta_1,\,\beta_2,\,\beta_3,\,0).$ This completes the proof.

\bigskip

\noindent {\bf Theorem 1:} If $a_i,\,i=1,\,2,\,3,\,4,$ are nonzero rational numbers,  the complete  solution of a solvable equation of type (\ref{qd4var}) can be expressed in terms of bilinear forms in four parameters if and only if 
$a_1a_2a_3a_4$ is a perfect square. If $a_1a_2a_3a_4=k^2$ and  $x_i=\alpha_i,\,i=1,\,2,\,3,\,4,$ is a nontrivial solution of equation (\ref{qd4var}) such that $a_3\alpha_3^2+a_4\alpha_4^2 \neq 0,$  this complete  bilinear   solution of  (\ref{qd4var}) 
 is given by 
\begin{equation}
\begin{array}{rcl}
x_1&=&(a_1a_3\alpha_1\alpha_3+k\alpha_2\alpha_4)a_1a_2a_4pr-(a_2a_3\alpha_2\alpha_3-k\alpha_1\alpha_4)a_1a_2a_4ps\\&& \;\;-(a_2a_3\alpha_2\alpha_3-k\alpha_1\alpha_4)a_1a_2a_4qr-(a_1a_3\alpha_1\alpha_3+k\alpha_2\alpha_4)a_2^2a_4qs,\\
x_2&=&(a_2a_3\alpha_2\alpha_3-k\alpha_1\alpha_4)a_1^2a_4pr+(a_1a_3\alpha_1\alpha_3+k\alpha_2\alpha_4)a_1a_2a_4ps\\ && \;\; +(a_1a_3\alpha_1\alpha_3+k\alpha_2\alpha_4)a_1a_2a_4qr-(a_2a_3\alpha_2\alpha_3-k\alpha_1\alpha_4)a_1a_2a_4qs,\\
x_3&=&-a_1a_2a_4(a_3\alpha_3^2+a_4\alpha_4^2)(a_1pr+a_2qs),\\
x_4&=&-a_1a_2k(a_3\alpha_3^2+a_4\alpha_4^2)(ps-qr), \label{bilin1}
\end{array}
\end{equation}
where $p,\,q,\,r,\,s$ are arbitrary parameters. 

Alternatively, if
$y_i=\beta_i,\,i=1,\,2,\,3,$ is a solution of equation (\ref{qd3var}) such that $\beta_3 \neq 0,$ the complete  bilinear   solution of  (\ref{qd4var}) 
 may be written  as follows:
\begin{equation}
\begin{array}{rcl}
x_1&=&a_4(a_1\beta_1p+a_2\beta_2q)r-a_2a_4(-\beta_2p+\beta_1q)s,\\
x_2 &=& a_1a_4(-\beta_2p+\beta_1q)r+a_4(a_1\beta_1p+a_2\beta_2q)s,\\
x_3&=&-a_4\beta_3(a_1pr+a_2qs),\\
x_4&=&k\beta_3(qr-ps) \label{bilin2}
\end{array}
\end{equation}
where, as before,  $p,\,q,\,r,\,s$ are arbitrary parameters.

\noindent {\bf Proof:}  It $a_1a_2a_3a_4=k^2, \; k \neq 0,$ and equation (\ref{qd4var}) is solvable, it follows from Lemma 1 that  equation (\ref{qd3var}) necessarily has a nontrivial solution  and we may take this solution as $y_i=\beta_i,\,i=1,\,2,\,3,$ with  $\beta_3 \neq 0.$  We will now  obtain  the complete solution of (\ref{qd4var}) by  using  the linear transformation given by
\begin{equation}
\begin{bmatrix} x_1 & x_2 &x_3 &x_4\end{bmatrix}=\begin{bmatrix} a_1a_4\beta_1& -a_2a_4\beta_1& a_2a_4\beta_2& a_2a_4\beta_2 \\ -a_1a_4\beta_2& a_2a_4\beta_2&a_1a_4\beta_1& a_1a_4\beta_1
 \\-a_1a_4\beta_3&-a_2a_4\beta_3&0&0\\ 0&0&-k\beta_3& k\beta_3 \end{bmatrix} \begin{bmatrix} X_1\\ X_2\\ X_3\\ X_4 \end{bmatrix}. \label{lttx}
 \end{equation}
 As the determinant of the matrix of this linear transformation  works out to $4ka_1a_2a_4^3\beta_3^2(a_1\beta_1^2+a_2\beta_2^2)$ which is not zero since $\beta_3 \neq 0$,  the linear transformation (\ref{lttx}) is invertible. Applying this transformation, equation (\ref{qd4var}) becomes 
 \begin{equation}
 \begin{array}{l}
a_4\{(a_1X_1-a_2X_2)^2+a_1a_2(X_3+X_4)^2)\}(a_1\beta_1^2+a_2\beta_2^2)\\
\quad +\{a_1^2a_3a_4X_1^2+2a_1a_2a_3a_4X_1X_2+a_2^2a_3a_4X_2^2\\
\quad \quad \quad \quad +k^2(X_3-X_4)^2\}\beta_3^2=0. \label{qd4varltt} \end{array}
\end{equation}
Using the relations $\beta_3^2=-(a_1\beta_1^2+a_2\beta_2^2)a_3$ (which is a consequence of $y_i=\beta_i,\,i=1,\,2,\,3,$ being a solution of (\ref{qd3var})), and $k^2=a_1a_2a_3a_4,$ equation (\ref{qd4varltt}) reduces to
\begin{equation}
-4a_1a_2a_4(a_1\beta_1^2+a_2\beta_2^2)(X_1X_2-X_3X_4) = 0,
\end{equation}
where we note that $4a_1a_2a_4(a_1\beta_1^2+a_2\beta_2^2) \neq 0.$ Thus, we get $X_1X_2-X_3X_4=0$ and as is well-known \cite[p. 69]{Si}, the complete solution of this equation  is given by
\begin{equation}
X_1=pr,\quad  X_2=qs, \quad X_3=ps, \quad X_4=qr. \label{valX} 
\end{equation}
 With these values of $X_i,$  the complete solution of (\ref{qd4var}) is given by (\ref{lttx}). This is a bilinear solution which is given  explicitly by (\ref{bilin2}). If $x_i=\alpha_i,\,i=1,\,2,\,3,$ is a nontrivial solution of equation (\ref{qd4var}) such that $a_3\alpha_3^2+a_4\alpha_4^2 \neq 0,$   we may take the values of $\beta_i,\,i=1,\,2,\,3,$ as defined by (\ref{valbeta}), and substituting these values in (\ref{bilin2}), we get the bilinear solution (\ref{bilin1}). 

We will now prove that  if equation (\ref{qd4var}) has a bilinear solution in four parameters, then  $a_1a_2a_3a_4$ must necessarily be a perfect square. If equation (\ref{qd4var}) has such a bilinear solution, it also has a linear solution in terms of two independent parameters, and therefore, without loss of generality,  we may take $x_1,\;x_2$ as arbitrary and  write this solution as $x_3=f_1x_1+f_2x_2,\;x_4=f_3x_1+f_4x_2$ where $f_i,\,i=1,\,2,\,3,\,4,$ are rational numbers and we note that both $f_1$ and $f_2$ cannot be 0. Substituting these values of $x_3$ and $x_4$ in (\ref{qd4var}), we get
 \begin{equation}
 (a_1+f_1^2a_3+f_3^2a_4)x_1^2+(2f_1f_2a_3+2f_3f_4a_4)x_1x_2+(a_2+f_2^2a_3+f_4^2a_4)x_2^2 = 0. \label{condx1x2}
 \end{equation}
 Since (\ref{condx1x2}) is identically true for all values of $x_1$ and $x_2,$ the coefficients of $x_1^2,\,x_1x_2$ and $x_2^2$ in (\ref{condx1x2}) must be 0, and we thus get three equations which can be solved for $a_1,\,a_2$ and $a_3.$ If $f_1$ and $f_2$ are both nonzero, the values of $a_1,\,a_2$ and $a_3$  lead to  the relation 
 \begin{equation}
 a_1a_2a_3a_4=a_4^4(f_1f_4-f_2f_3)^2f_3^2f_4^2/(f_1^2f_2^2). \label{condbilin}
 \end{equation}
 If $f_1=0$ then for equation (\ref{condx1x2}) to be identically satisfied for all $x_1,\,x_2,$ we must necessarily also have $f_4=0$ and  we then get $a_1=-f_3^2a_4,\,a_2=-f_2^2a_3,$  so that $a_1a_2a_3a_4=f_2^2f_3^2a_3^2a_4^2.$ Similarly, if $f_2=0,$ we get $f_3=0,\,a_1=-f_1^2a_3,\,a_2=-f_4^2a_4$ and so  $a_1a_2a_3a_4=f_1^2f_4^2a_3^2a_4^2.$ Thus in all cases, $a_1a_2a_3a_4$ is a perfect square. This completes the proof.
 \bigskip

\noindent  {\bf  Theorem 2:} If $Q(x_1,\,x_2,\,x_3,\,x_4)$ is a quadratic form, with rational coefficients, in four independent variables $x_1,\,x_2,\,x_3,\,x_4,$ the matrix of the form $Q(x_i)$ is $A$ and the diophantine equation
 \begin{equation}
 Q(x_1,\,x_2,\,x_3,\,x_4)=0, \label{qdgen4var}
 \end{equation}
 has a solution in integers, the complete solution of (\ref{qdgen4var}) can be expressed in terms of bilinear forms in four variables if and only if the determinant  $|A|$ is a nonzero perfect square.
 
\noindent  {\bf Proof:} It is well-known that any rational quadratic form in $n$ variables can be reduced by an invertible rational linear transformation to a diagonal form. Denoting the column vectors $\{x_1,\,x_2,\,x_3,\,x_4\}$ and $\{y_1,\,y_2,\,y_3,\,y_4\}$ by $X$ and $Y$ respectively,  the transpose of any matrix $M$ by $M^{\prime}$, and the invertible linear tranformation that reduces the quadratic form $Q(x_1,\,x_2,\,x_3,\,x_4)$ to a diagonal form by $X=PY$,  we have $|P| \neq 0,$ and 
 \begin{equation}
 P^{\prime}AP=\mbox{diag} \begin{bmatrix}a_1 & a_2 &a_3 &a_4\end{bmatrix}, \label{quadred}
 \end{equation}
where $a_i \in {\mathbb Q}\setminus\{0\}$ (since the form $Q(x_1,\,x_2,\,x_3,\,x_4)$  has four independent variables), while  equation (\ref{qdgen4var}) reduces to the equation
\begin{equation}
a_1y_1^2+a_2y_2^2+a_3y_3^2+a_4y_4^2=0.
\label{qd4vary}
\end{equation}
  It follows from (\ref{quadred}) that   $ | P^{\prime}|. | A | .  |P| = a_1a_2a_3a_4, $  and so,
  \begin{equation} 
   |P|^2. |A|=a_1a_2a_3a_4. \label{relPAai}
   \end{equation} 
  
  If the diophantine equation (\ref{qdgen4var}) is solvable, its complete solution can be expressed in terms of bilinear forms in four variables if and only if the complete  solution of   equation (\ref{qd4vary}) can be so expressed. By Theorem 1, this  is possible if and only if $a_1a_2a_3a_4$ is a nonzero  perfect square, and the theorem now follows readily from (\ref{relPAai}).
	
	A bilinear solution of equation (\ref{qdgen4var}) may be obtained by reducing it, using a suitable linear transformation $X=PY$, to an equation of type (\ref{quadred}), obtaining a bilinear solution of (\ref{quadred}) by Theorem 1, and finally using the relations $X=PY$.

  We note en passant that it  follows from (\ref{relPAai}) that $|A|=0$ if and only if the form $Q(x_1,\,x_2,\,x_3,\,x_4)$ has fewer than four independent variables. When $|A|=0,$  at least one of the $a_i$ must be 0. If only one of the four rational numbers $a_i$ is  0,
  then equation (\ref{qd4vary}) and so also equation (\ref{qdgen4var}), cannot have a linear solution in two independent parameters. If two of these four numbers, say $a_3$ and $a_4$ are 0 and $-a_1a_2$ is a perfect square, it is readily seen that both equation (\ref{qd4vary}) and  equation (\ref{qdgen4var}) have a linear solution in three independent parameters, and the same is true if three of the four numbers $a_i$ are 0.
	
	We now prove a theorem which shows that even under the weaker assumption that equation (\ref{qdgen4var}) has a solution in linear forms in two variables, it is true that $|A|$ is a nonzero perfect square.
\bigskip

\noindent  {\bf Theorem 3:} If $Q(x_1,\,x_2,\,x_3,\,x_4)$ is a quadratic form, with rational coeffcients, in four independent variables $x_1,\,x_2,\,x_3,\,x_4,$ the matrix of the form $Q(x_i)$ is $A$ and the diophantine equation (\ref{qdgen4var}) 
 has a solution in terms of linear forms in two independent variables,  then $|A|$ is a nonzero perfect square and the complete solution of (\ref{qdgen4var}) can be expressed in terms of bilinear forms in four variables. 
 
 \noindent{\bf Proof:} As in Theorem 2, we apply an invertible  linear transformation $X=PY$ to reduce equation (\ref{qdgen4var}) to the diagonal form 
(\ref{qd4vary}), and further, we  have $|P|^2. |A|=a_1a_2a_3a_4.$ Equation (\ref{qd4vary})  will also have a solution in terms of linear forms in two independent variables, and since we have already shown   in the last part the proof of Theorem 1 that in such a situation,  $a_1a_2a_3a_4$ must be a nonzero perfect square, it follows that $|A| =a_1a_2a_3a_4|P|^{-2}$  must also be a nonzero perfect square. It now follows from  Theorem 2  that the complete solution of (\ref{qdgen4var}) can be expressed in terms of bilinear forms in four variables.

\subsection{Two numerical examples}

\hspace{0.25in} We now give two numerical examples illustrating the application of the results of Section 2.1. As a simple example, we consider the quadratic equation

\begin{equation}
x_1^2-9x_2^2-x_3^2+4x_4^2=0, \label{ex1}
\end{equation}
 which has a solution $(x_1,\,x_2,\,x_3,\,x_4)=(1,\,0,\,1,\,0),$ and the product $a_1a_2a_3a_4$ is a perfect square, namely 36, so the conditions of Theorem 1 are satisfied.  Using (\ref{bilin1}) and replacing $r,\,s$ respectively by $m/18$ and $n/18$ respectively, we obtain the following complete solution of (\ref{ex1}): 
\begin{equation}
\begin{array}{rclrcl}
x_1 &=& 2pm+18qn, &x_2 &=& 2pn+2qm, \\x_3 &=& -2pm+18qn, &x_4 &=& -3pn+3qm, \label{solex1}
\end{array} \end{equation}
where $p,\,q,\,m,\,n$ are arbitary parameters.

As a second example, we consider the equation (\ref{ex2})  
which has a solution $(x_1,\,x_2,\,x_3,\,x_4)=(0,\,2,\,-1,\,-1),$ and the determinant of its  matrix $A$ works out to 36 which is a perfect square. Applying the linear transformation,
\begin{equation}
\begin{array}{rclrcl}
x_1 &=& y_1-y_2+y_3,& x_2 &=& 31y_1-29y_2+26y_3,\\ x_3 &=& -19y_1+18y_2-16y_3, &x_4 &=& y_4, \label{lttex2}
\end{array}
\end{equation}
reduces equation (\ref{ex2}) to the diagonal form,
\begin{equation}
y_1^2+2y_2^2-6y_3^2-3y_4^2=0. \label{ex2y}
\end{equation}
Equation (\ref{ex2y}) has a solution $(y_1,\,y_2,\,y_3,\,y_4)=(1,\,1,\,0,\,1)$, and applying Theorem 1, we get a bilinear solution of (\ref{ex2y}) which, on substituting $r=-m/18,s=-n/36,$ may be written as follows:
\begin{equation}
\begin{array}{rclrcl}
y_1 &=& (2p+2q)m+(p-2q)n, &y_2 &=& (-p+2q)m+(p+q)n, \\y_3 &=& pm+qn, &y_4 &=& -pn+2qm, \label{solyex2}
\end{array}
\end{equation}
where $p,\,q,\,m,\, n,$ are arbitrary parameters. Now using the relations (\ref{lttex2}), we get the bilinear solution  (\ref{solex2}) of the equation (\ref{ex2}).

We note here that the complete solution of equation (\ref{ex2}), obtained by the usual method  as described in \cite[p. 432]{Di}, is as follows:
\begin{equation}
\begin{array}{rcl}
x_1&=&16pq+26pr,
 \\x_2 &=& 540p^2+152pq+304pr-16q^2-70qr-70r^2, \\
x_3 &=& -270p^2-76pq-152pr+16q^2+64qr+61r^2,\\
x_4 &=& -270p^2-76pq-152pr+16q^2+48qr+35r^2, \label{solex2b}
\end{array}
\end{equation}
where $p,\,q,\,r$ are arbitrary parameters. The  bilinear solution (\ref{solex2}) is also a complete solution and  is clearly much simpler than the complete solution (\ref{solex2b}).

\section{A necessary condition for the solvability of a pair of quadratic equations}

\hspace{0.25in}In this section we consider the solvability of a pair of simultaneous quadratic diophantine equations in four variables. Naturally each of the  two equations must be individually solvable for otherwise the question of their having a common solution does not arise. We first prove a theorem that gives a necessary, though not sufficient, condition for the solvability of  such  equations  and then  we describe a method of obtaining all integer solutions when this condition is satisfied. 

\noindent  {\bf Theorem 4:} A necessary condition that the simultaneous equations  in four independent variables 
 \begin{eqnarray}
   Q_1(x_1,\,x_2,\,x_3,\,x_4) &=&0, \label{qd1}\\
    Q_2(x_1,\,x_2,\,x_3,\,x_4) &=&0,\label{qd2}
    \end{eqnarray}
    have a nontrivial solution in integers is that there exists a nontrivial solution in integers of the  quartic equation,
    \begin{equation}
    \eta^2=f(\xi_1,\,\xi_2)= |A(\xi_1,\,\xi_2)|, \label{eqnxi12}
    \end{equation}
    where $A(\xi_1,\,\xi_2)$ is the matrix of the quadratic form $\xi_1Q_1(x_1,\,x_2,\,x_3,\,x_4) +\xi_2Q_2(x_1,\,x_2,\,x_3,\,x_4).$ Further, a necessary condition that there exist at least two distinct solutions of the simultaneous equations (\ref{qd1}) and (\ref{qd2}) is that there exists an integer  solution of equation (\ref{eqnxi12}) with $\eta \neq 0.$
    
    \noindent {\bf Proof:} Any rational solution $ (\xi_1,\,\xi_2,\,\eta)$ of (\ref{eqnxi12}) implies the existence of an integer  solution $ (k\xi_1,\,k\xi_2,\,k^2\eta)$ for a suitably chosen value of $k$, and hence wherever required,  it suffices to prove just the existence of a rational solution of (\ref{eqnxi12}).

		We will first show that when the  simultaneous quadratic equations (\ref{qd1}) and (\ref{qd2}) have two distinct integer solutions, then (\ref{eqnxi12}) has a solution in integers with $\eta \neq 0.$ Let $A_j, \,j=1,\,2,$ be the respective matrices of the quadratic forms $ Q_j(x_i),\,j=1,\,2.$  If either of the determinants $|A_i|, $ say $|A_1|$ is  a nonzero perfect square, say  $\eta_1^2$,  a rational   solution of equation (\ref{eqnxi12}) is given by $(\xi_1,\,\xi_2,\,\eta)=(1,\,0,\,\eta_1).$  
    
    Next we assume that   neither  $|A_1|$ nor  $|A_2|$ is a  perfect square.  If  $(x_1,\,x_2,\,x_3,\,x_4) = (\alpha_1,\,\alpha_2,\,\alpha_3,\,\alpha_4)$ and  $(x_1,\,x_2,\,x_3,\,x_4) = (\beta_1,\,\beta_2,\,\beta_3,\,\beta_4)$ be  two integer solutions of equations  (\ref{qd1}) and (\ref{qd2}), we substitute 
    \begin{equation}
    x_i=\alpha_im+\beta_in, \quad i=1,\,2, \,3,\,4, \label{subsx}
    \end{equation}
 in the equation
 \begin{equation}
 \xi_1Q_1(x_1,\,x_2,\,x_3,\,x_4) +\xi_2Q_2(x_1,\,x_2,\,x_3,\,x_4) =0, \label{eqnxiQ12}
 \end{equation}
 when we get the equation
 \begin{equation}
 mn(\xi_1h_1+\xi_2h_2)=0, \label{eqnxih}
 \end{equation}
 where both $h_1$ and $h_2$ are  determined by $\alpha_i,\,i=1,\,2, \,3,\,4,$ and $\beta_i,\,i=1,\,2, \,3,\,4. $ We note that both $h_1$ and $h_2$ must be nonzero since either $h_j$  being 0 would imply that (\ref{subsx}) is a linear solution in two independent parameters $m$ and $n$  of the corresponding equation  $Q_j(x_1,\,x_2,\,x_3,\,x_4)=0,$  and in view of Theorem 3, the corresponding determinant $|A_i|$ must be a nonzero perfect square contradicting our assumption.  Taking $\xi_1=-h_2$ and $\xi_2=h_1$, we observe that a  linear solution of equation (\ref{eqnxiQ12}) in terms of two independent parameters $m$ and $n$ is given by (\ref{subsx}) and hence it follows from  Theorem 3 that with these values of $\xi_j,\,j=1,\,2,$ the determinant $A(\xi_1,\,\xi_2)$ must be a nonzero perfect square, that is, the diophantine equation (\ref{eqnxi12}) has a nontrivial solution with $\eta \neq 0.$

		 We will now consider the case when   the simultaneous equations (\ref{qd1}) and (\ref{qd2}) have  just  one  solution in integers, say $(x_1,\,x_2,\,x_3,\,x_4) = (\alpha_1,\,\alpha_2,\,\alpha_3,\,\alpha_4)$ where we assume without loss of generality that $\alpha_4 \neq 0.$  By the  invertible linear transformation defined  by $x_1=y_1+\alpha_1y_4,\,x_2=y_2+\alpha_2y_4,\,x_3=y_3+\alpha_3y_4,\,x_4=\alpha_4y_4,$  we transform these two equations to the corresponding equations $Q^{\prime}_j(y_1,\,y_2,\,y_3,\,y_4)=0,\,j=1,\,2$ which have just one solution, namely $(y_1,\,y_2,\,y_3,\,y_4)=(0,\,0,\,0,\,1)$, and hence the two quadratic forms $Q^{\prime}_j(y_i),\,j=1,\,2$ may be written as follows:
		\begin{equation}
		\begin{array}{rcl}
		Q^{\prime}_1(y_1,\,y_2,\,y_3,\,y_4)& = & y_4L_1(y_1,\,y_2,\,y_3)+Q_3(y_1,\,y_2,\,y_3),\\
		Q^{\prime}_2(y_1,\,y_2,\,y_3,\,y_4)& = & y_4L_2(y_1,\,y_2,\,y_3)+Q_4(y_1,\,y_2,\,y_3),
		\end{array}
		\end{equation}
		where $L_j(y_1,\,y_2,\,y_3),\;j=1,\,2$ are linear forms and $Q_j(y_1,\,y_2,\,y_3),\;j=3,\,4$ are quadratic forms respectively in the variables $y_1,\,y_2,\,y_3.$ There are now two possibilities:
		
		\medskip
		
		\noindent (i) If the two linear forms $L_j(y_1,\,y_2,\,y_3),\;j=1,\,2$ are linearly dependent, there exist integers $h_1,\,h_2$, both not zero,  such that $h_1L_1(y_1,\,y_2,\,y_3)+h_2L_2(y_1,\,y_2,\,y_3)=0,$ and so the quadratic form $h_1	Q^{\prime}_1(y_1,\,y_2,\,y_3,\,y_4)+h_2Q^{\prime}_2(y_1,\,y_2,\,y_3,\,y_4)$ does not contain $y_4$ and therefore has at most three independent variables.  Thus the determinant of the matrix of this quadratic form must be $0$, and hence the determinant of the matrix of the quadratic form $h_1Q_1(x_1,\,x_2,\,x_3,\,x_4)+h_2Q_2(x_1,\,x_2,\,x_3,\,x_4)$ is also $0$. Thus, $ |A(h_1,\,h_2)|=0$, and $(\xi_1,\,\xi_2,\,\eta)=(h_1,\,h_2,\,0) $ is a nontrivial solution of the equation (\ref{eqnxi12}). 
		
		\medskip
		
		\noindent (ii) If the two linear forms $L_j(y_1,\,y_2,\,y_3),\;j=1,\,2$ are linearly independent, we choose a suitable arbitrary  third form $L_3(y_1,\,y_2,\,y_3)$ such that the three forms $L_j(y_1,\,y_2,\,y_3),\;j=1,\,2,\,3$ are linearly independent, and apply the invertible linear transformation defined by $L_j(y_1,\,y_2,\,y_3)=z_j,\;j=1,\,2,\,3,$  and $y_4=z_4$ to the two quadratic forms $Q^{\prime}_j(y_i),\;j=1,\,2,$ to get corresponding quadratic forms $Q^{\prime \prime }_j(z_i)\;j=1,\,2,$ which may be written as follows:
		\begin{equation}
		\begin{array}{rcl}
		Q^{\prime \prime}_1(z_i)& = & z_4z_1+a_1z_1^2+a_2z_2^2+a_3z_3^2+a_4z_2z_3+a_5z_1z_3+a_6z_1z_2,\\
		Q^{\prime \prime}_2(z_i)& = & z_4z_2+b_1z_1^2+b_2z_2^2+b_3z_3^2+b_4z_2z_3+b_5z_1z_3+b_6z_1z_2.
		\end{array}
		\end{equation}	
where the coefficients $a_j,\,b_j,\; j=1,\,\ldots,\,6,$ can be effectively computed. The two equations $Q^{\prime \prime}_j(z_i)=0,\,j=1,\,2$ necessarily have just one solution in integers, namely,  $(z_1,\,z_2,\,z_3,\,z_4)=(0,\,0,\,0,\,1)$ and hence both $a_3$ and $b_3$ cannot be $0.$ Next we consider the equation $b_3Q^{\prime \prime}_1(z_i)-a_3Q^{\prime \prime}_2(z_i)=0.$ On substituting $z_2=z_1$ and cancelling out the factor $z_1,$ this equation reduces to a linear equation  and hence the equation $b_3Q^{\prime \prime}_1(z_i)-a_3Q^{\prime \prime}_2(z_i)=0$ has a  linear solution in two independent parameters.  It follows that the corresponding equation $b_3Q_1(x_i)-a_3Q_2(x_i)=0$ also has a linear solution in two independent parameters, and hence by Theorem 3, the determinant of the matrix of the quadratic form $b_3Q_1(x_i)-a_3Q_2(x_i)$ must be a nonzero perfect square. Thus, $|A(b_3,\,-a_3)|$ is a nonzero perfect square and  hence equation (\ref{eqnxi12}) has a nontrivial solution. 

Thus even when the simultaneous equations (\ref{qd1}) and (\ref{qd2}) have a single solution in integers, the diophantine  equation (\ref{eqnxi12}) has a nontrivial  solution in integers. This completes the proof of the theorem.

It is interesting to note that  the simultaneous equations (\ref{qd1}) and (\ref{qd2}) may essentially have just one numerical solution while the corresponding equation (\ref{eqnxi12}) has infinitely many solutions. As an example, the two equations,
\begin{equation}
\begin{array}{rcl}
 Q_1(x_1,\,x_2,\,x_3,\,x_4) &= &x_1x_2-x_3x_4=0,\\
Q_2(x_1,\,x_2,\,x_3,\,x_4)&=&(x_1-x_2)^2+(x_1-x_3)^2+(x_1-x_4
)^2=0,
\end{array}
\end{equation}
have only one solution, namely $(x_1,\,x_2,\,x_3,\,x_4)=(1,\,1,\,1,\,1)$ while  the determinant $A(\xi_1,\,\xi_2)|=\xi_1^2(\xi_1+2\xi_2)(\xi_1-6\xi_2)/16,$ and hence equation (\ref{eqnxi12}) is easily seen to have infinitely many solutions. 

We have noted earlier that the condition of Theorem 4 is not sufficient for the solvability of the two equations. This is illustrated by the two equations,
\begin{equation}
\begin{array}{rcl}
 Q_1(x_1,\,x_2,\,x_3,\,x_4) &= &x_1x_2-x_3x_4=0,\\
Q_2(x_1,\,x_2,\,x_3,\,x_4)&=&(x_1-x_2)^2+(x_1-x_3)^2 +(7x_1-x_4)^2=0,
\end{array}
\end{equation}
which clearly do not have a nontrivial solution in integers even though both are individually solvable,  the determinant $|A(\xi_1,\,\xi_2)|=(\xi_1^4-4\xi_1^3\xi_2-204\xi_1^2\xi_2^2-96\xi_1\xi_2^3)/16,$ and the equation (\ref{eqnxi12}) has infinitely many solutions, two of them being $(\xi_1,\,\xi_2,\,\eta)=(0,\,1,\,0)$ and $(\xi_1,\,\xi_2,\,\eta)=(54,\,2,\,549).$
  
  \subsection{A method of solving  a pair of simultaneous diophantine equations in four  variables}
 
   \hspace{0.25in} We will now use Theorem 4 to determine the solvability of a given   pair of simultaneous diophantine equations (\ref{qd1}) and (\ref{qd2}), and  show that  the complete solution of these equations may be obtained by solving the equivalent system consisting of either (\ref{qd1}) or (\ref{qd2})  and the equation 
    \begin{equation}
    \xi_1Q_1(x_1,\,x_2,\,x_3,\,x_4)+\xi_2 Q_2(x_1,\,x_2,\,x_3,\,x_4)=0 \label{qd12}
    \end{equation}
    for appropriately chosen $\xi_1$ and $\xi_2.$ 
		
		As before, let $A(\xi_1,\,\xi_2)$ be the matrix of the quadratic form $\xi_1Q_1(x_i)+\xi_2 Q_2(x_i).$ We now consider the solvability of the equation 
		(\ref{eqnxi12}). If it has no solutions in integers, it follows from Theorem 4 that the pair of equations (\ref{qd1}) and (\ref{qd2}) has no solutions. As an example, when 
		\begin{equation}
		\begin{array}{rcl}
		Q_1(x_i)&=& x_1^2+2x_2^2-x_3^2-x_4^2,\\
		Q_2(x_i)&=&7x_1^2+4x_1x_2+14x_2^2-6x_3^2+2x_3x_4-8x_4^2,
		\end{array}
		\end{equation}
 both  the  equations $Q_j(x_i)=0,\;j=1,\,2,$ individually have infinitely many solutions, but $|A(\xi_1,\,\xi_2)|=32(\xi_1^2+14\xi_1\xi_2+47\xi_2^2)^2$ and so equation (\ref{eqnxi12}) does not have any nontrivial solutions. Hence these simultaneous equations are not solvable. 
		
		When equation (\ref{eqnxi12}) is solvable, we will choose a suitable integer solution of (\ref{eqnxi12}), solve the corresponding equation (\ref{qd12}), and use this solution to obtain the complete solution of equations (\ref{qd1}) and (\ref{qd2}). During this process, we may need to solve quartic equations of the type 
    \begin{equation}
    \eta^2=c_0\xi^4+c_1\xi^3+c_2\xi^2+c_3\xi+c_4, \label{quarticeta}
    \end{equation}
    where $c_i,\,i=1,\,\ldots,\,4,$ are rational numbers.   While there is no algorithm that will give a numerical solution of an equation of type (\ref{quarticeta}) in every case, such a  solution can often be found using programs such as APECS.
		
		The following two subsections,   based on the reducibility of the quartic form $|A(\xi_1,\,\xi_2)|$, discuss the method further and give  illustrative  examples.

  \subsubsection{The determinant $|A(\xi_1,\,\xi_2)|$ has a linear factor} 
	\hspace{0.25in}In this case it is simplest to  choose integer values of $\xi_1,\,\xi_2$ such that $|A(\xi_1,\,\xi_2)|=0$ so that  equation (\ref{qd12}) has fewer than four independent variables. We now apply  an invertible linear transformation $X=PY$ to  reduce this equation to the diagonal form (\ref{qd4vary}) where one or more of the coefficients $a_i$ must be zero. We obtain the complete solution of this equation 
	(including its trivial solutions!), determine the values of $x_i,\; i=1,\,\ldots,\, 4,$  by using the relations $X=PY$, substitute these values in $(\ref{qd1})$ (or in $(\ref{qd2})$) and solve the resulting equation  to obtain the complete solution of the simultaneous equations (\ref{qd1}) and (\ref{qd2}). 
	
	There are several possibilities depending on the solvability of the reduced equation (\ref{qd4vary}). We give four examples  illustrating different types of solutions that can be obtained.  

	For the simultaneous equations,
	\begin{equation}
	\begin{array}{rclcl}
	Q_1(x_i)&= &x_1^2+5x_2^2-4x_2x_4-3x_3^2+2x_4^2&=&0,\\
  Q_2(x_i)&= &x_1^2+3x_2^2-2x_2x_4-2x_3^2+x_4^2&=&0, \label{qdpairex1}
  \end{array}
	\end{equation}
we have $|A(\xi_1,\,\xi_2)|=-(\xi_1+\xi_2)(2\xi_1+\xi_2)(3\xi_1+2\xi_2)^2,\; 2Q_1(x_i)-3Q_2(x_i)=-(x_1+x_2-x_4)(x_1-x_2+x_4),$ and the complete solution is given by
	\[
  x_1 = (m^2-n^2),\; x_2 = 2mn,\; x_3 = (m^2+n^2),\; x_4 = (m^2+2 m n-n^2),\] and
   \[x_1 = (m^2-n^2),\; x_2 = 2mn,\; x_3 = (m^2+n^2),\; x_4 =- (m^2-2mn-n^2),\]
	where $m$ and $n$ are arbitrary parameters.

	As a second example, for the simultaneous equations,
	\begin{equation}
	\begin{array}{rclcl}
	Q_1(x_i)&= &x_1^2+544x_1x_2-320x_1x_3-27x_2^2&& \\&& \quad \quad \quad +320x_2x_4-x_3^2+320x_4^2&=&0,\\
  Q_2(x_i)&= & x_1^2+1088x_1x_2-640x_1x_3-55x_2^2&& \\&& \quad \quad \quad +640x_2x_4-x_3^2+640x_4^2&=&0, \label{qdpairex2}
  \end{array}
	\end{equation}
	we have $|A(\xi_1,\,\xi_2)|=320(\xi_1+2\xi_2)(299\xi_1+599\xi_2)(97\xi_1+193\xi_2)^2,\; 2Q_1(x_i)-Q_2(x_i)=x_1^2+x_2^2-x_3^2,$ and the complete solution is given by the single numerical solution $(x_1,\,x_2,\,x_3,\,x_4)=(3,\,4,\,5,\,-2),$ and the two parametric solutions
  \[\begin{array}{rcl}
  x_1 &=& -119m^4-480m^3n-350m^2n^2+25n^4,\\ x_2 &=& 120m^4+200m^3n+120m^2n^2+200mn^3,\\ x_3 &=& 169m^4+480m^3n+450m^2n^2+25n^4, \\x_4 &=& -155m^4-186m^3n+300m^2n^2+70mn^3-25n^4,
  \end{array}\]and 
   \[\begin{array}{rcl}
  x_1 &=& -119m^4-480m^3n-350m^2n^2+25n^4,\\ x_2 &=& 120m^4+200m^3n+120m^2n^2+200mn^3,\\ x_3 &=& 169m^4+480m^3n+450m^2n^2+25n^4,\\ x_4 &=& 35m^4-14m^3n-420m^2n^2-270mn^3+25n^4,
   \end{array}\]
	where, as before,  $m$ and $n$ are arbitrary parameters.

As a third example, for the simultaneous equations,
	\begin{equation}
	\begin{array}{rclcl}
	Q_1(x_i)&= & x_1^2+4x_2^2+8x_2x_3+8x_2x_4&& \\&&\quad \quad +5x_3^2+16x_3x_4+8x_4^2&=&0,\\
  Q_1(x_i)&= &2x_1^2+5x_2^2+8x_2x_3+8x_2x_4&& \\&&\quad \quad +4x_3^2+16x_3x_4+8x_4^2&=&0, \label{qdpairex3}
  \end{array}
	\end{equation}
	we 
	obtain just four solutions, namely  $(\pm2,\,0,\,2,\,-1)$ and $(\pm2,\,0,\,2,\,-3).$

Finally,  the  pair of  quadratic equations,
  \begin{equation}\begin{array}{rclcl}
  Q_1(x_i)&=& x_1^2+6x_2^2+2x_2x_3+16x_2x_4 &&\\&&\quad \quad -4x_3^2+8x_3x_4+16x_4^2&=&0,\\
  Q_2(x_i)&= &2x_1^2+7x_2^2+2x_2x_3+16x_2x_4&& \\&&\quad \quad-5x_3^2+8x_3x_4+16x_4^2&=&0, \label{qdpairex4}
  \end{array}\end{equation}
  have infinitely many solutions that are given by
  \[ x_1 = 4(\xi^2-1),\; x_2 = 8\xi,\; x_3 = 4(\xi^2+1),\; x_4 = -\xi^2-4\xi \pm 2\eta -1,\]
  where $\xi$ and $\eta$ are related by
  \begin{equation} \eta^2=\xi^4+\xi^3+\xi^2+\xi+1.\label{ec1etaxi}\end{equation}
  Equation (\ref{ec1etaxi}) represents an elliptic curve of rank 1 with a point of infinite order being given by $(\xi,\,\eta)=(-1,\,1).$ We can thus find infinitely many solutions of (\ref{ec1etaxi}) and hence also of the above pair of quadratic equations.

 \subsubsection{The determinant $|A(\xi_1,\,\xi_2)|$ does not have a linear factor} 
\hspace{0.25in}In this case we take $\xi_2=\xi \xi_1$ and replace $\eta$ by $\eta \xi_1^2$ when (\ref{eqnxi12}) reduces to an equation of type (\ref{quarticeta}) and we have to determine its solvability.  If (\ref{eqnxi12}) has no integer solutions, Theorem 4 ensures that the simultaneous equations (\ref{qd1}) and (\ref{qd2}) have no integer solutions. If (\ref{eqnxi12}) has an integer solution, we choose integer values of $\xi_1,\,\xi_2$ such that $|A(\xi_1,\,\xi_2)|$ is a perfect square, and  obtain the complete solution of  equation (\ref{qd12}) in terms of bilinear forms. On substituting this solution in  (\ref{qd1}), we get  a parametrised quadratic equation  whose solvability depends on the solvability of another  quartic equation of type (\ref{quarticeta}). We thus determine the solvability of the  simultaneous equations (\ref{qd1}) and (\ref{qd2}) and obtain their complete solution.

  As an example, we consider the simultaneous equations,
   \begin{equation}
	\begin{array}{rclcl}
  Q_1(x_i)&= &x_1^2-9x_2^2-x_3^2+4x_4^2&=&0, \\
   Q_2(x_i)&= &3x_1^2-30x_1x_2-4x_1x_3-9x_2^2\\&& -12x_2x_3-7x_3^2+12x_3x_4+4x_4^2&=&0. \label{qdpairex5}
   \end{array} \end{equation}
	This is the pair of equations (\ref{qdpairex50}) mentioned in the Introduction. Here $|A(\xi_1,\,\xi_2)|= 36(\xi_1^2+6\xi_1\xi_2+20\xi_2^2)^2$ and so any linear combination of the two equations has a bilinear solution. We note that the equation $Q_1(x_i)=0$ is the same as equation (\ref{ex1}), and its solution is given by (\ref{solex1}). 
On substituting the values of $x_i$ given by (\ref{solex1}) in $Q_2(x_i)=0$ and solving the resulting equation, we obtain the complete solution of  equations (\ref{qdpairex5}) as given by (\ref{qdpairex5sol1}) and (\ref{qdpairex5sol2}).

 As a second example,  we consider the two simultaneous equations,
 \begin{equation} \begin{array}{rclcl}
Q_1(x_i)&=& 6x_1^2-4x_1x_2+4x_2x_3-36x_2x_4-5x_3^2-27x_4^2&=&0, \\
Q_2(x_i)&=& 11x_1^2-8x_1x_2+4x_2^2+8x_2x_3-72x_2x_4-9x_3^2&& \\ && \quad \quad \quad \quad\quad \quad\quad \quad \quad\quad \quad\quad \quad\quad\quad \quad -63x_4^2&=&0. \label{qdpairex6}
 \end{array}
\end{equation}
Here we get $|A(\xi_1,\,\xi_2)|= 9828\xi_1^4+78336\xi_1^3\xi_2+234360\xi_1^2\xi_2^2+311616\xi_1\xi_2^3+155268\xi_2^4$ and taking $(\xi_1,\,\xi_2)=(2,\,-1)$ makes $|A(\xi_1,\,\xi_2)|$ a perfect square so that the equation $2Q_1(x_i)-Q_2(x_i)=0$ has a bilinear solution which is given by
\[x_1 = 6pm+6qn, x_2 = -3pm+3qn, x_3 = 6qm+6pn, x_4 = 2qm-2pn.\]
 On substituting these values of $x_i$  in either of the two equations (\ref{qdpairex6}), we get the quadratic equation
\[(2p^2+pq-2q^2)m^2+(-2p^2+2pq-q^2)mn+(-2p^2+2pq+q^2)n^2=0,\]
which has a rational solution if and only if $p,\,q$ are chosen such that the discriminant $20p^4-16p^3q-24p^2q^2+8pq^3+9q^4$ becomes a perfect square. We thus obtain the following solution of the simultaneous equations (\ref{qdpairex6}):
 \[\begin{array}{rcl}
 x_1 &=& 12\xi^3+12\xi^2+18\xi+6\eta\xi-24,\\ x_2 &=& -6\xi^3+18\xi^2+3\xi-3\eta\xi-12,\\ x_3 &=& 24\xi^3+24\xi^2-36\xi+6\eta+6,\\ x_4 &=& -8\xi^3+4\xi+2\eta+2,
 \end{array}\] 
 where $\xi$ and $\eta$ are related by $\eta^2=20\xi^4-16\xi^3-24\xi^2+8\xi+9.$ This last quartic equation represents an elliptic curve of rank 3, three points of infinite order being $(\xi,\,\eta)=(1/2,\,5/2),$ $(\xi,\,\eta)=(-1/2,\,3/2),$ and $(\xi,\,\eta)=(-14/95,\,4897/1805).$  We thus get infinitely many solutions of  equations  (\ref{qdpairex6}), the solution arising from the first of the aforementioned three points being $(x_1,\,x_2,\,x_3,\,x_4)=(6,\,21,\,-24,\,-16).$

 As a final example, 
 we consider the  simultaneous equations,
 \begin{equation} \begin{array}{rclcl}
Q_1(x_i)&=&x_1^2-2x_1x_2-9x_2^2+3x_3^2-4x_3x_4+11x_4^2&=&0, \\
Q_2(x_i)&=&6x_1^2-x_1x_2+x_2^2-15x_3^2-2x_3x_4-11x_4^2 &=&0, \label{qdpairex7}
 \end{array} \end{equation}
both of which individually have infinitely many solutions. Here we get $|A(\xi_1,\,\xi_2)|= -290\xi_1^4+454\xi_1^3\xi_2+(37739/4)\xi_1^2\xi_2^2-(20035/2)\xi_1\xi_2^3+943\xi_2^4$, and taking 
$\xi_1=4,\xi_2=3$ makes $|A(\xi_1,\,\xi_2)|$ a perfect square so that the equation $4Q_1(x_i)+3Q_2(x_i)=0$ has a bilinear solution  which is given by 
\[x_1 = 12pm+4qn, x_2 = 8pm-4qn, x_3 = 5qm+5pn, x_4 = -5qm+15pn.\]
On substituting these values of $x_i$  in  either of the two equations (\ref{qdpairex7}), we get the quadratic equation
\[(75q^2-104p^2)m^2-166pqmn+(375p^2-16q^2)n^2=0,\]
which has a rational solution if and only if $p,\,q$ are chosen such that the discriminant $400(390p^4-229p^2q^2+12q^4)$ becomes a perfect square. It is readily established by elementary congruence considerations that this discriminant can never become a perfect square for integer values of $p$ and $q$. Thus the  simultaneous  equations  (\ref{qdpairex7})  do not have any integer solutions.

\section{Applications to  quartic diophantine equations}

\hspace{0.25in}We will now use bilinear solutions of quadratic diophantine equations to obtain  rational solutions of the quartic 
diophantine equation
\begin{equation}
y^2=t^4+a_1t^3+a_2t^2+a_3t+a_4, \label{quartic1}
\end{equation}
where the coefficients $a_i,\,i=1,\,2,\,3,\,4,$ are rational numbers. This equation represents a quartic model of an elliptic curve and has been studied extensively. Various methods of finding rational   solutions of this equation are discussed in \cite[pp. 639-644]{Di} and \cite[pp. 69-70, 77-78]{Mo}. However,  explicit formulae  have not been published for deriving new rational solutions  of (\ref{quartic1}) starting from  one or two known solutions of this equation. The reason for this seems to be that the existing methods lead to cumbersome formulae.    

In Section 4.1 we prove a theorem that gives formulae by which new rational solutions of (\ref{quartic1}) may be obtained starting from either a single rational solution of (\ref{quartic1}) or from two such solutions.  In Section 4.2 we will show  these  formulae generate rational solutions of (\ref{quartic1}) that are much simpler than the solutions obtained by existing methods. Finally in Section 4.3 we show that when two solutions of the more general quartic equation
\begin{equation}
y^2=a_0t^4+a_1t^3+a_2t^2+a_3t+a_4, \label{quartic2}
\end{equation}
are known, we can still use  the theorem  of Section 4.1   to find new rational solutions of (\ref{quartic2}). As in the case of equation (\ref{quartic1}), these solutions are also much simpler than the ones obtained by existing methods.

\subsection{Rational solutions of the equation $y^2=t^4+a_1t^3+a_2t^2+a_3t+a_4$}

\noindent {\bf Theorem 5:} If $(t_1,\,y_1)$ is any rational solution of equation (\ref{quartic1}) with $y_1 \neq 0$, a new rational solution of  (\ref{quartic1}) is given by $(t_{11},\,y_{11})$ where
\begin{equation}
\begin{array}{rl}
t_{11}=&\{-4(2t_1^2+2a_1t_1+a_2)y_1^2+4(2t_1^4+a_1t_1^3-a_3t_1-2a_4)y_1\\& \quad \quad \quad  +(4t_1^3+3a_1t_1^2+2a_2t_1+a_3)^2\} \\ & \times [4y_1\{(4t_1+a_1)y_1+4t_1^3+3a_1t_1^2+2a_2t_1+a_3\}]^{-1},\\
y_{11}=&\{64y_1^6+(128t_1^2+64a_1t_1-16a_1^2+64a_2)y_1^5+(320t_1^4+320a_1t_1^3\\ & \;\;+(96a_1^2+64a_2)t_1^2+64(a_1a_2-a_3)t_1-16a_1a_3+16a_2^2)y_1^4\\& \;\;
+8(4t_1^3+3a_1t_1^2+2a_2t_1+a_3)(16t_1^3+12a_1t_1^2+3a_1^2t_1\\ &\;\;+a_1a_2-2a_3)y_1^3-2(4t_1+a_1)(4t_1^3+3a_1t_1^2+2a_2t_1+a_3)^3y_1\\ & \;\;-(4t_1^3+3a_1t_1^2+2a_2t_1+a_3)^4 \}
\\& \quad \times [4y_1\{(4t_1+a_1)y_1+4t_1^3+3a_1t_1^2+2a_2t_1+a_3\}]^{-2}, \label{t11}
\end{array}
\end{equation}
provided that $\{(4t_1+a_1)y_1+4t_1^3+3a_1t_1^2+2a_2t_1+a_3\} \neq 0.$

Further, if   $(t_1,\,y_1)$ and $(t_2,\,y_2)$ are two rational solutions of equation (\ref{quartic1}) with $t_1 \neq t_2$, a new rational solution of  (\ref{quartic1}) is given by $(t_{12},\,y_{12})$ where
\begin{equation}
\begin{array}{rl} 
t_{12}=&\{-2y_1y_2+2(t_1-t_2)(t_2y_1-t_1y_2)+a_1(t_1+t_2)t_1t_2\\& \quad  +2a_2t_1t_2+a_3(t_1+t_2)+2a_4+2(t_1^2-t_1t_2+t_2^2)t_1t_2\}\\& \quad \quad \times \{(t_1-t_2)(2y_1-2y_2+a_1(t_1-t_2)+2t_1^2-2t_2^2)\}^{-1},\\
y_{12}=&[-(y_1-y_2)^4-(t_1-t_2)\{(4t_1+a_1)y_1-(4t_2+a_1)y_2\}(y_1-y_2)^2\\&  -(t_1-t_2)^3\{(2t_1+2t_2+a_1)(y_1^2-y_2^2)-4(t_1-t_2)y_1y_2\} \\ & 
+(t_1-t_2)^3\{(4t_1+a_1)(t_1^2+2t_1t_2+3t_2^2+a_1t_1+2a_1t_2+a_2)y_1\\& \quad \quad  -(4t_2+a_1)(3t_1^2+2t_1t_2+t_2^2+2a_1t_1+a_1t_2+a_2)y_2\}\\& 
+(t_1-t_2)^4(t_1^2+2t_1t_2+3t_2^2+a_1t_1+2a_1t_2+a_2)\\ &\quad  \times (3t_1^2+2t_1t_2+t_2^2+2a_1t_1+a_1t_2+a_2)]\\ & \quad \quad \times \{(t_1-t_2)(2y_1-2y_2+a_1(t_1-t_2)+2t_1^2-2t_2^2)\}^{-2}, \label{t12}
\end{array}
\end{equation}
provided that $\{2y_1-2y_2+a_1(t_1-t_2)+2t_1^2-2t_2^2\} \neq 0.$

\noindent {\bf Proof:} If    $(t_1,\,y_1)$ is a solution of (\ref{quartic1}), we have $a_4=y_1^2-(t_1^4+a_1t_1^3+a_2t_1^2+a_3t_1)$, and with this value of $a_4$, it is readily verified using a symbolic computation  software such as MAPLE that $(t_{11},\,y_{11})$ is a solution of (\ref{quartic1}). Similarly, if $(t_1,\,y_1)$ and $(t_2,\,y_2)$ are two  solutions of equation (\ref{quartic1}), we solve the relations $y_i^2=t_i^4+a_1t_i^3+a_2t_i^2+a_3t_i+a_4,\; i=1,\,2$ for $a_3,\;a_4$, and using these values of $a_3,\;a_4$, it is readily verified by  MAPLE that $(t_{12},\,y_{12})$ is a solution of (\ref{quartic1}).

While direct computation suffices to prove the theorem, the following more elaborate proof gives a genesis of the formulae (\ref{t11}) and (\ref{t12}). 

We will first derive the formulae (\ref{t12}) that give a new solution of (\ref{quartic1}) starting from two known rational solutions of this equation. 

We  make an initial assumption that the coefficients $a_i$ satisfy the condition $a_1^3-4a_1a_2+4a_3 \neq 0.$ We will now establish a relationship between   the rational solutions of equation (\ref{quartic1}) and  the  solutions of the simultaneous quadratic diophantine equations given below:
 \begin{eqnarray}
Q_1(x_i)= x_1^2+c_1x_2^2+c_2x_2x_3+x_2x_4+c_3x_3x_4+c_4x_4^2&=&0,\label{qd5}\\
 Q_2(x_i)=x_1x_2-x_3x_4&=&0, \label{qd6}
 \end{eqnarray}
 where the values of $c_j,\;j=1,\,2,\,3,\,4,$ are given below:
 \begin{equation}
 \begin{array}{rcl}
 c_1 &=& (a_1^2-4a_2)/16,\\ c_2 &=& -(a_1^3-4a_1a_2+4a_3)/16,\\ c_3 &=& -a_1/2, \\c_4 &=&4(a_1^4-4a_1^2a_2+8a_1a_3-16a_4)/(a_1^3-4a_1a_2+4a_3)^2. \label{defcj}
 \end{array}
 \end{equation}

Let $(\alpha_1,\,\alpha_2,\,\alpha_3,\,\alpha_4)$ be an arbitrary rational solution of the simultaneous  equations (\ref{qd5}) and (\ref{qd6})  and let $[\alpha_1,\,\alpha_2,\,\alpha_3,\,\alpha_4]$, written briefly as $\bar{\alpha},$ denote the equivalence class of all rational solutions $(k\alpha_1,\,k\alpha_2,\,k\alpha_3,\,k\alpha_4)$ where $ k \in {\mathbb Q} \setminus\{0 \}.$  Let $S$ be the set of all  equivalence classes  $[\alpha_1,\,\alpha_2,\,\alpha_3,\,\alpha_4]$ where $\alpha_4 \neq 0$ and  let $T$ be the set of all rational solutions $(t_1,\,y_1)$ of equation (\ref{quartic1}). We will show that there is a one-one correspondence between the sets $S$ and $T$.

Let $\bar{\alpha}$ be any arbitrary element of $S$. We define a  mapping $\phi: S \rightarrow T$  as follows: 
\[\phi(\bar{\alpha})= (t_0,\,y_0)\]
where
\begin{equation}
\begin{array}{rcl}
t_0&=&\frac{\displaystyle (a_1^3-4a_1a_2+4a_3)\alpha_2+8a_1\alpha_4}{\displaystyle 16\alpha_4}, \\[0.2in]
y_0&=&\frac{\displaystyle (a_1^3-4a_1a_2+4a_3)\{32\alpha_1\alpha_4-(a_1^3-4a_1a_2+4a_3)\alpha_2^2-8a_1\alpha_2\alpha_4\}}{\displaystyle 256\alpha_4^2}. \label{defphi}
\end{array}
\end{equation}

We note that the mapping $\phi $  is well-defined since it is easily seen that the values of $t_0$ and $y_0$ are independent of the scalar factor $k$. We will now show that $(t_0,\,y_0)$ is indeed a solution of equation (\ref{quartic1}).  Since  $\bar{\alpha} \in S, $ we have $Q_1(\alpha_i)=0$ as well as  $Q_2(\alpha_i)=0,$ and on eliminating $\alpha_3$ from these two relations, we find that $\alpha_1,\,\alpha_2,\,\alpha_4$ satisfy the condition $f(\alpha_1,\,\alpha_2,\,\alpha_4)=0$ where 
\[
\begin{array}{l}
 f(\alpha_1,\,\alpha_2,\,\alpha_4)=(a_1^3-4a_1a_2+4a_3)^3\alpha_1\alpha_2^2-16(a_1^3-4a_1a_2+4a_3)^2\alpha_2\alpha_4^2\\\quad \quad \quad\quad\quad \quad\quad\;-(64a_1^4-256a_1^2a_2+512a_1a_3-1024a_4)\alpha_4^3\\\quad \quad \quad\quad\quad \quad\quad\;-(a_1^3-4a_1a_2+4a_3)^2(\alpha_2^2a_1^2-8\alpha_1\alpha_2a_1-4\alpha_2^2a_2+16\alpha_1^2)\alpha_4.
\end{array}
\]
Now, using the values of $(t_0,\,y_0)$ given by (\ref{defphi}), we find that 
\[y_0^2-(t_0^4+a_1t_0^3+a_2t_0^2+a_3t_0+a_4)=f(\alpha_1,\,\alpha_2,\,\alpha_4)/(1024\alpha_4^3)=0.\]
This shows that $(t_0,\,y_0) $ is a solution of equation (\ref{quartic1}) and hence $ \phi(\bar{\alpha}) \in T.$

Next let $(t_1,\,y_1)$ be any arbitrary element of $T$. We define a  mapping $\psi: T \rightarrow S$  as follows:
\[\psi(t_1,\,y_1)=[\beta_1,\,\beta_2,\,\beta_3,\,\beta_4],\]
where
\begin{equation}
\begin{array}{rcl}
\beta_1&=& 4(a_1^3-4a_1a_2+4a_3)(2t_1^2+t_1a_1+2y_1),\\
\beta_2&=&-8(a_1^3-4a_1a_2+4a_3)(2t_1+a_1),\\
\beta_3&=&-32(2t_1^2+t_1a_1+2y_1)(2t_1+a_1),\\
\beta_4&=&(a_1^3-4a_1a_2+4a_3)^2. \label{defpsi}
\end{array}
\end{equation}
We will now show that $[\beta_1,\,\beta_2,\,\beta_3,\,\beta_4] \in S.$ It is readily verified that $Q_1(\beta_i)=0.$ Further, with the above values of $\beta_i$, we get
\[Q_2(\beta_i)=64(a_1^3-4a_1a_2+4a_3)^2\{y_1^2-(t_1^4+a_1t_1^3+a_2t_1^2+a_3t_1+a_4)\}=0.\]
Thus, $(\beta_1,\,\beta_2,\,\beta_3,\,\beta_4)$ is  a solution of the simultaneous equations (\ref{qd5}) and (\ref{qd6}). Finally we note that, in view of our initial ssumption,  $\beta_4 \neq 0.$  This proves that $[\beta_1,\,\beta_2,\,\beta_3,\,\beta_4] \in S.$

With the above definitions of the mappings $\phi$ and $\psi$, it is readily established that for any  $(t_1,\,y_1) \in T,$ we have $\phi(\psi(t_1,\,y_1))=(t_1,\,y_1)$ and for any 
$\bar{\alpha} \in S,$ we have $\psi(\phi(\bar{\alpha}))=\bar{\alpha}.$ We have thus established a one-to-one correspondence between the two sets $S$ and $T$.

Now, let  $(t_1,\,y_1)$ and  $(t_2,\,y_2)$ be any two distinct  elements of $T$ such that $t_1 \neq t_2$. Further, let $\psi(t_1,\,y_1)=[\beta_1,\,\beta_2,\,\beta_3,\,\beta_4] $ and $\psi(t_2,\,y_2) =[\gamma_1,\,\gamma_2,\,\gamma_3,\,\gamma_4],$ so that $(\beta_1,\,\beta_2,\,\beta_3,\,\beta_4)$ and $(\gamma_1,\,\gamma_2,\,\gamma_3,\,\gamma_4)$ are two distinct solutions of the simultaneous equations (\ref{qd5}) and (\ref{qd6}). Thus, we have $Q_j(\beta_1,\,\beta_2,\,\beta_3,\,\beta_4)=0, \; j=1,\,2$ as well as $Q_j(\gamma_1,\,\gamma_2,\,\gamma_3,\,\gamma_4)=0, \; j=1,\,2.$ The values of $\beta_i,\,i=1,\,2,\,3,\,4,$ are given by (\ref{defpsi}) while the values of  $\gamma_i,\,i=1,\,2,\,3,\,4,$ are given by analogous relations with $t_1,\,y_1$ being replaced by $t_2,\,y_2$ respectively. 

We now consider the quadratic equation 
\begin{equation}
Q_1(x_1,\,x_2,\,x_3,\,x_4) +tQ_2(x_1,\,x_2,\,x_3,\,x_4)=0, \label{qd56} 
\end{equation}
where $t$ is some rational number. It is readily verified that the determinant $|A|$ of the matrix of the quadratic form $Q_1(x_i) +tQ_2(x_i)$ is $t^4+a_1t^3+a_2t^2+a_3t+a_4.$ We will now use the solutions $(\beta_1,\,\beta_2,\,\beta_3,\,\beta_4)$ and $(\gamma_1,\,\gamma_2,\,\gamma_3,\,\gamma_4)$ of equations (\ref{qd5}) and (\ref{qd6}) to find a value of $t$ for which equation (\ref{qd56}) has a solution in terms of linear forms in two independent variables.  By Theorem 3, this value of $t$ will make the determinant $|A|$ a perfect square, and we will thus obtain a solution of equation (\ref{quartic1}). 

 We make the following substitution in  equation (\ref{qd56}):
    \begin{equation}
		\begin{array}{rclrcl}
    x_1 &=& \beta_1m+\gamma_1n,\quad&  x_2 &=& \beta_2m+\gamma_2n,\\ x_3 &=& \beta_3m+\gamma_3n, \quad & x_4 &=& \beta_4m+\gamma_4n. \label{subsxx}
    \end{array}  \end{equation}
   Since $Q_j(\beta_i)=0, \; j=1,\,2$  and also $Q_j(\gamma_i)=0, \; j=1,\,2,$ equation (\ref{qd56}) reduces to 
\begin{equation}
\begin{array}{l}
\{(\beta_1\gamma_2+\beta_2\gamma_1-\beta_3\gamma_4-\beta_4\gamma_3)t+2\beta_1\gamma_1+(2c_1\beta_2+c_2\beta_3+\beta_4)\gamma_2\\
\quad \quad +(c_2\beta_2+c_3\beta_4)\gamma_3+(\beta_2+c_3\beta_3+2c_4\beta_4)\gamma_4\}mn=0. \label{eqnt} 
\end{array}
\end{equation}		

We now choose a value of $t$ such that equation (\ref{eqnt}) is satisfied. We find, using the values of 	 $\beta_i,\gamma_i$ and the values of the coefficients $c_j, \;j=1,\,2,\,3,\,4,$ given by (\ref{defcj}), that it is possible to choose such a value of $t$ when $\{2y_1-2y_2+a_1(t_1-t_2)+2t_1^2-2t_2^2\} \neq 0 $ and that this value of $t$ is given by $t_{12}$  as defined by (\ref{t12}). Thus, when $t=t_{12}$, a solution of (\ref{qd56}) is given by (\ref{subsxx}) in terms of two independent linear parameters $m$ and $n$. It follows from Theorem 3  that  the determinant of   the matrix  of the quadratic form $Q_1(x_i) +tQ_2(x_i)$, that  is,  $t^4+a_1t^3+a_2t^2+a_3t+a_4$ must be a perfect square. Thus, $t=t_{12}$ gives a rational solution of (\ref{quartic1}). 

Using the value of $t_{12}$ obtained above, the corresponding value of $y_{12}$ given by (\ref{t12}) was obtained by  direct computation using MAPLE. 

We now consider the quartic equation (\ref{quartic1}) when the coefficients $a_i$ are such that $a_1^3-4a_1a_2+4a_3 = 0$ so that our initial assumption does not hold. In this case, we make a change of variable by replacing $t$ by $t+m$  where the  value of $m$ is suitably chosen such that coefficients of the  transformed quartic equation do not satisfy the condition $a_1^3-4a_1a_2+4a_3 = 0. $ We then use the known solutions $(t_1,\,y_1)$ and $(t_2,\,y_2)$ of equation (\ref{quartic1}) to find corresponding known solutions of the transformed quartic equation, and use them by applying the formula just found   to find a  new rational solution of the transformed quartic equation, and eventually obtain a new solution of the original quartic equation (\ref{quartic1}). This new  solution is given by  $(t_{12}, \, y_{12})$ so that  the formulae (\ref{t12}) hold even when  $a_1^3-4a_1a_2+4a_3 = 0.$ The only condition is that $(t_1-t_2)\{2y_1-2y_2+a_1(t_1-t_2)+2t_1^2-2t_2^2\} \neq 0.$

To obtain formulae (\ref{t11}), we consider the value of $t_{12}$ given by (\ref{t12}) as a continuous function of $t_2$ while $t_1$ and $y_1$ are fixed. We found that 
$ \mbox {\rm lim}_{t_2 \rightarrow t_1}t_{12} = t_{11},$
where the value of $t_{11}$ is given by (\ref{t11}). Direct computation now confirmed that $t=t_{11}$  indeed gives a rational solution of (\ref{quartic1}) and also yielded the value of $y_{11}$ as given by (\ref{t11}). This completes the proof.
 
\noindent {\bf Corollary 1:}  The quartic function of $t$ given by \[(t-a_1)(t-a_2)(t^2+a_3t+a_4) \] becomes a perfect square if we take 
$ t=  (a_1a_2-4)/(a_1+a_2+a_3)$. 
This follows immediately from (\ref{t12}) by taking $(t_1,\,y_1)=(a_1,\,0)$ and $(t_2,\,y_2)=(a_2,\,0)$ as the two known solutions.

\noindent {\bf Corollary 2:}  The quartic function of $t$ given by \[(t-a_1)(t-a_2)(t^2+a_3t+a_4)+m^2 \] becomes a perfect square if we take  
$ t=   (a_1a_2-2m-4)/(a_1+a_2+a_3) $. 
As in the case of Corollary 1 , this follows immediately from (\ref{t12}) by taking $(t_1,\,y_1)=(a_1,\,m)$ and $(t_2,\,y_2)=(a_2,\,m)$ as the two known solutions.

 We note that the formulae (\ref{t12}) are symmetric with respect to the two solutions $(t_1,\,y_1)$ and $(t_2,\,y_2)$ and  so interchanging these solutions   in formulae (\ref{t12}) does not yield a new   solution. Further, since the existence of any solution $(t_1,\,y_1)$ of equation (\ref{quartic1}) automatically implies the existence of the solution $(t_1,\,-y_1)$ , we can choose either  sign for  $y_1$ and $y_2$ in (\ref{t11}) and (\ref{t12}) and in general, obtain two  new solutions of (\ref{quartic1}) starting from a single known solution and four new solutions starting from two known solutions. There are, of course, examples when we get fewer  solutions, or in exceptional situations, no new solutions. 

\subsection{A comparison with existing methods}

\hspace{0.25in}We will now compare the results of Theorem 5 with the existing methods of obtaining rational solutions of  equation (\ref{quartic1}). 

 Fermat gave a procedure (as quoted by Dickson \cite[p. 639]{Di}) of finding a rational solution of (\ref{quartic1}) under fairly general conditions as well as 
of finding a new rational solution of  (\ref{quartic2}) starting from a known solution. The procedure can often, though not always, be repeated any number of times to obtain infinitely many  solutions of (\ref{quartic1}). The solution of (\ref{quartic1}) obtained by Fermat's method is given  by a rational function of the coefficients $a_i$ with the numerator being a polynomial of degree 4 in $a_i$ and the denominator being a polynomial of degree 3  in $a_i$.  Therefore  repetitions of the process  rapidly generate rational solutions of great height.  

Mordell \cite[p. 70]{Mo} gives a substitution based on two known rational solutions $(t_1,\,y_1)$ and $(t_2,\,y_2)$ of  (\ref{quartic2}) and this  yields  several new solutions of  (\ref{quartic2}).   While the resulting new rational solutions are not given explicitly, they can be worked out using a software such as MAPLE and are given by  rational functions of $(t_1,\,y_1)$ and $(t_2,\,y_2)$. For any given quartic equation of type (\ref{quartic1}), if we take $t_2=O(t_1)$, Mordell's method gives one set of  solutions in which  the numerator is of  the order of $t_1^9$ while the denominator is of order of $t_1^8$, and   another  set of solutions   in which  the numerator is of   $O(t_1^{17})$ while the denominator is of  $O(t_1^{16})$. 

We can also apply a birational transformation to the quartic equation (\ref{quartic1}) to transform it to a cubic equation (see \cite[p. 77]{Mo}), obtain two rational solutions of this cubic equation corresponding to  two known solutions of (\ref{quartic1}), use these two solutions to obtain another rational solution of the cubic equation and finally find the corresponding solution of (\ref{quartic1}) by using the birational transformation. This leads to very complicated formulae where again assuming  that $t_2=O(t_1)$, we find  that  the numerator of the new solution is of  the order of $t_1^9$ while the denominator is of order of $t_1^8$. 

As compared to the above,  for any given quartic equation of type (\ref{quartic1}), if we take $t_2=O(t_1)$,   the  numerator of the solution (\ref{t12}) given by Theorem 5 is of order of $t_1^4$ while the denominator is of order of  $t_1^3. $

 We now consider  a numerical example. The quartic equation
\begin{equation}
y^2 = t^4-7t^3-3t^2+48t-35, \label{quartex1}
\end{equation}
has the  solutions $(t,\,y)=(1,\,\pm 2)$ and $(t,\,y)=(2,\,\pm 3)$.   On applying formulae (\ref{t12}) to these solutions, we get four values of $t$  that make the right-hand side of (\ref{quartex1}) a perfect square, namely, $18,\;-30/11,\;58/9$ and $-22/3,$ while applying  formulae (\ref{t11}) to the known  solutions individually, we get four additional such values of $t,$ namely, 
$1001/152,\;-729/248,\;-121/39$ and $421/57.$ 

As a comparison,  we now find the values of $t$ that make the right-hand side of (\ref{quartex1}) a perfect square by various other methods. Direct application of Fermat's method gives $t=-5961/344$ while using the known solutions and applying other classical methods as discussed in \cite[pp. 639-640]{Di}, we obtain the following values of $t$: $221/97,\;141/65,\;142/109,$ and $70/61.$ Repetition of the procedures yields solutions  of greater height. 

The substitution  mentioned by Mordell leads to the following values of $t$: $1002/425,$ $ 2202/2009,\;169/112,\; 1201/520,\;1547585/829473$ and $52238/61093.$

For  a further comparison, we used the software APECS to determine rational solutions of (\ref{quartic1}). First, we reduced the quartic equation (\ref{quartic1}), by a birational tansformation,  to the Weierstrass minimal form,  
\begin{equation}
Y^2 = X^3-199X-29. \label{weier}
\end{equation} 
Equation (\ref{weier}) represents an elliptic curve of rank 4, a basis for this elliptic curve being given by the four points $P_1=(-1, 13),\;P_2=(-5, 29),\; P_3= (-13, 19),\; P_4=(15, 19).$ The four values of $t$ corresponding to these four points are as follows: $70/61,\;6/7,\;142/109$ and $18,$  
 while the four values of $t$ corresponding to the points $2P_i,\;i=1,\,2,\,3,\,4,$ are as follows: $877941/130585,\;[11365/377, \;-2703129/708263, \;-21785/4807.$  
Adding two of the four points $P_i$ at a time, we get six new points on the elliptic curve (\ref{weier}), and the values of $t$ corresponding to these six points are as follows: $223/31,\;-5961/344,$
$ 7401/4112,\;-2287/680,\;55/24,\;141/65.$ 

On observing all the solutions obtained by various methods, it  seems to be quite clear  that the solutions obtained by the formulae of Theorem 5  compare very well with the existing methods. The comparison is even better  when the coefficients $a_i$ of the quartic function are given in parametric terms, and one or two solutions of (\ref{quartic1}) are known, for then the formulae of Theorem 5  easily  give new solutions of lower degree as compared to other methods. 

\subsection{Rational solutions of the equation $y^2=a_0t^4+a_1t^3+a_2t^2+a_3t+a_4$} 
\hspace{0.25in} We now show that is possible to use Theorem 5 even in the case of the more general quartic equation (\ref{quartic2}) provided we know two solutions $(t_1,\,y_1)$ and $(t_2,\,y_2)$ of (\ref{quartic2}) such that $t_1 \neq t_2.$  In fact, on substituting $t=s+t_1,$ in (\ref{quartic2}), we get $y^2=b_0s^4+b_1s^3+b_2s^2+b_3s+y_1^2,$ where the coefficients $b_j$ are determined by the coefficients $a_i$ and $t_1$, and now the further substitutions $y=y_1Y/r^2,\; s=1/r,$ yield an equation of the type $Y^2=r^4+c_1r^3+c_2r^2+c_3r+c_4$ for which we can readily find a solution corresponding to the known solution  $(t_2,\,y_2)$ of (\ref{quartic2}). We  have thus  transformed equation (\ref{quartic2}) to an equation of type (\ref{quartic1}) with a known solution, and can now  apply the theorem. We illustrate the method by an example.

Consider the equation,
\begin{equation}
y^2=2t^4+3t^3+7t^2-207t+379, \label{quartex2}
\end{equation}
for which two known solutions are $(2,\,7)$ and $(3,\,8).$ The substitution $t=s+2$ reduces the equation to  
\begin{equation}
y^2 = 2s^4+19s^3+73s^2-79s+49,\label{quartex2a}
\end{equation}
and the further substitutions $y=7Y/(r^2),s=1/r,$  yield the quartic equation
\begin{equation}
Y^2= r^4-(79/49)r^3+(73/49)r^2+(19/49)r+2/49. \label{quartex2b}
\end{equation}
A solution of the quartic equation (\ref{quartex2b}), corresponding to the known solution $(3,\,8)$ of equation (\ref{quartex2}) is $(1,\,8/7)$, and now using formula (\ref{t11}), we get two  solutions of (\ref{quartex2b}), leading to two solutions of (\ref{quartex2}) given by   $t=122/17$ and $t=1754/809.$ 

Following the classical method, after obtaining equation (\ref{quartex2a}), we would write $y=7-(79/14)s+(8067/2744)s^2,$ so that the constant term and the coefficients of $s$ and $s^2$ on both sides become equal, and this gives a solution of  (\ref{quartex2a}) with $s=9582440/1219937,$ and we finally get a solution of  (\ref{quartex2a}) given by  $t=12022314/1219937.$ 

This example shows that even in the case of the more general equation (\ref{quartic2}), applying Theorem 5 yields simpler solutions as compared to the classical methods

\noindent Postal Address: 13/4 A Clay Square, \\
\hspace{1.1in} Lucknow - 226001,\\
\hspace{1.1in}  INDIA \\
\medskip
\noindent e-mail address: ajaic203@yahoo.com


\begin{thebibliography}{9}

\bibitem{Di}L. E. Dickson, {\it History of theory of numbers, }   Vol. 2,
Chelsea Publishing Company, New York, 1992, reprint.
\bibitem{Mo} L. J. Mordell, {\it Diophantine equations}, Academic Press, London (1969). 
 \bibitem{Si} W. Sierpinski (edited by A. Schinzel), {\it Elementary Theory of Numbers}, PWN-Polish Scientific Publishers, Warszawa (1987). 


 \end{thebibliography}
\end{document}